\documentclass{article}
\usepackage{amsfonts,amssymb, amsmath}

\newtheorem{prop}{Proposition}

\def\bq{ \begin{equation}}
\def\eq{ \end{equation}}
\def\ben{ \begin{eqnarray}}
\def\en{ \end{eqnarray}}

\begin{document}


\title{Multisymplectic structure of nonintegrable H\'{e}non-Heiles system}
\author{A. V. Tsiganov,\\
International Laboratory for Mirror Symmetry \\
National Research University Higher School of Economics \\
Moscow, Russia\\
e--mail: andrey.tsiganov@gmail.com}
\date{}
\maketitle

\begin{abstract}
Multi-symplectic integrators are typically regarded as a discretization of the Hamiltonian partial differential equations. 
This is due to the fact that, for generic finite-dimensional Hamiltonian systems, there exists only one independent symplectic structure. In this note, the second invariant symplectic form is presented for the nonintegrable H\'{e}non-Heiles system, Kepler problem, integrable and non-integrable Toda type systems. This approach facilitates the construction of a multi-symplectic integrator, which effectively preserves both symplectic forms for these benchmark problems.
 \end{abstract}

\section{Introduction}
In the majority of cases of dynamical systems, it is not possible to find an analytical solution. This increases the necessity for advanced numerical techniques to facilitate model-based design and analysis. As computing power increases, it becomes possible to produce numerical solutions over longer time intervals. It is essential to ensure that the qualitative properties of the integrator are fully understood in order to guarantee the accuracy of the numerical simulation and the reliability of long-range predictions. Traditional numerical schemes do not account explicitly for the fundamental features of the underlying dynamical system, however, incurring error that may suggest non-physical behavior. Structure-preserving integrators are numerical methods that respect the fundamental physics of a problem by preserving the geometric properties of the governing differential equations \cite{mar2001,book0,book1,book2}.

The preservation of invariant structures of differential equations is a key factor in the effectiveness of the structure-preserving integrator. The discrete Lagrangian flow obtained by a classical variational integrator preserves a symplectic form. From this property it follows, by backward error analysis, that the energy is approximately preserved. Multisymplectic variational integrators are structure-preserving numerical schemes which preserve exactly the momenta associated with the symmetries, it is symplectic in time, and the energy is well conserved \cite{br01,rat16,sv2003}. 
For the St\"{a}ckel system there are structure-preserving integrators conserving the same number of constants of motion
as the degree of freedom \cite{min06}.

So, it can be posited that the preservation of the original system of differential equations' more invariant structures during discretization is indicative of the corresponding integrator's enhanced efficiency and stability. For instance, for the Kepler problem symplectic integrators provide good long-time behavior of the solution because they preserve the symplectic structure,  an approximate Hamiltonian function and the angular momentum \cite{book0,book1}. However, these structure-preserving integrators neither preserve Runge-Lenz vector even approximately.  Thus, development of novel discretisation methods has been undertaken for the purpose of preserving all invariants: energy, angular momentum, Runge-Lenz vector and symplectic structure, see \cite{kepl18}  and references therein.

In this note we consider a vector field $X$ on the manifold $M$ with coordinates $x=(x^1,\ldots,x^m)$  which defines a
system of ordinary differential equations 
\bq\label{m-eq}
\frac{d}{dt}\,x = X(x^1,\ldots,x^m)\,.
\eq 
Prior to the development of a structure-preserving integrator for (\ref{m-eq}), it is necessary to define the underlying
 geometric structures $T$ that have to be preserved. These underlying invariant structures are solutions of the  invariance equation
\bq\label{g-inv}
\mathcal L_X\, T=0
\eq
on the tensors invariants  $T$ to the flow generated by $X$,  including  phase space functions  (first
integrals), multivector fields (symmetry fields, Poisson structures),  differential forms (symplectic form, volume form),  etc. Here $L_X  T$ is a Lie derivative of tensor field $T$ along the vector field $X$.

The study of invariant tensor fields is a well-established area of research at least for the integrable Hamiltonian systems, see \cite{bog96,bog97,bog98, koz19,rat24}.  We will solve the invariance equation (\ref{g-inv})  for some known integrable and non-integrable models using a brute force method. Firstly, the H\'{e}non-Heiles system is selected  because it is a well-studied benchmark in the field of nonlinear dynamics \cite{hh64,hh25}.
Secondly, the Kepler problem is chosen as an example of degenerate or superintegrable systems.
Thirdly, the Toda type systems are chosen  as an example of Hamiltonian systems with non-polynomial potentials. 

\subsection{Hamiltonian systems}
The Lie derivative $L_X T$ determines the rate of change of the tensor field $T$ under the phase space deformation defined by the flow of the system (\ref{m-eq}). In local coordinates the Lie derivative of the tensor field $T$ of type $(p, q)$ is equal to
\begin{align*}
({\mathcal {L}}_{X}T)^{i_{1}\ldots i_{p}}_{j_{1}\ldots j_{q}}=\sum_{k=1}^n X^{k}(\partial _{k}T^{i_{1}\ldots i_{p}}_{j_{1}\ldots j_{q}})
&-\sum_{\ell=1}^n (\partial _{\ell}X^{i_{1}})T^{\ell i_{2}\ldots i_{p}}_{j_{1}\ldots j_{s}}-\ldots - \sum_{\ell=1}^n (\partial _{\ell}X^{i_{p}})T^{i_{1}\ldots i_{p-1}\ell}_{j_{1}\ldots j_{s}}\label{lie-d}\\
&+\sum_{m=1}^n(\partial _{j_{1}}X^{m})T^{i_{1}\ldots i_{p}}_{m j_{2}\ldots j_{q}}+\ldots +\sum_{m=1}^n(\partial _{j_{q}}X^{m})T^{i_{1}\ldots i_{p}}_{j_{1}\ldots j_{q-1}m}
\nonumber
\end{align*}
where $\partial _s= {\partial }/{\partial x^s}$ is the partial derivative on the $x^s$ coordinate.

The Lie derivative commutes with the exterior differentiation operation and satisfies the Leibnitz rule. It allows us to construct tensor invariants from a set of basic invariant tensor fields which either have a simpler functional dependence on the variables $x$, or have some special properties or physical interpretation. As an example, the Hamiltonian vector field $X$ on a $2n$-dimensional symplectic manifold is defined by
\[
\iota_X\omega=dH
\]
where $\iota$ is an interior product, $\omega$ is a symplectic 2-form and $dH$ is a 1-form constructed by differentiating the Hamilton function $H$ for which we suppose that 
\[\mathcal L_XH=0\qquad \mbox{and}\qquad \mathcal L_X\omega=0\,.\]
It means that $H$ is the scalar solution of (\ref{g-inv}), which usually coincides with the mechanical energy of the dynamical system (\ref{m-eq}), and $dH$ and $\omega$ are tensor solutions of the invariance equation (\ref{g-inv})  \cite{mbook}. 

Symplectic structure induces multisymplectic structures \cite{leon99,bar18}. So,using tensor product of the basic invariants $H$ and $\omega$ we can define a family of invariant differential forms of type $(0,2k)$ and $(0,2k-1)$
\bq\label{w-inv}
\omega_{2k}=\omega^k \qquad\mbox{and}\qquad \omega_{2k-1}=\iota_X\omega_{2k}\,,\qquad k=1,\ldots n\,,
\eq
and invariant multivector fields $P_j=\omega_j^{-1}$ of type $(j,0)$. If the discretization scheme preserves basic invariants $\omega \equiv\omega_2$ and $H$, it also preserve their respective derivative invariants (\ref{w-inv}).

Similar for the Hamiltonian vector fields on the Poisson manifold
\[
X=PdH\,,
\]
so that 
\[\mathcal L_XH=0\qquad \mbox{and}\qquad \mathcal L_X P=0\,,\]
where $P$ is an invariant Poisson bivector, we can construct a family of 
invariant multivector fields of type $(2k,0)$ and $(2k-1,0)$
\bq\label{p-inw}
P_{2k}=P^k\qquad\mbox{and}\qquad P_{2k-1}=P_{2k}dH.
\eq
Here $P^k$ is a tensor product of $k$ copies of bivector $P\equiv P_2$, which could be degenerate \cite{per96}.

The following general question arises in connection with the above: do systems
of Hamiltonian differential equations (\ref{m-eq}) admit non-trivial tensor invariants of a type $(p, q)$ which can not be obtained from the basic invariants?

We aim  to provide a direct response to this question by solving an invariance equation (\ref{g-inv}) 
for a classical non-relativistic system with two degrees of freedom on the plane by using a modern mathematical software. Let us take Hamiltonian and canonical Poisson bivector   
\bq\label{p-can}
H=\frac{1}{2}(p_1^2+p_2^2)+V(q_1,q_2)\,,
\qquad
P=\left(
 \begin{array}{cccc}
 0 & 0 & 1 & 0 \\
 0 & 0 & 0 & 1 \\
 -1 & 0 & 0 & 0 \\
 0 & -1 & 0 & 0
 \end{array}
 \right)\,,
\eq 
which are invariants to flow of the Hamiltonian vector field
\bq\label{x-vec}
X=PdH=\left(
 \begin{array}{c}
 p_1 \\
 p_2 \\ 
 -{\partial_1 V} \\ 
 -{\partial_2 V}
 \end{array}
 \right)\,.
\eq
In coordinates $x=(q_1,q_2,p_1,p_2)$ the invariance equation $\mathcal L_X P'=0$ (\ref{g-inv}) has the  form
\[
(L_X P')^{ij}=\sum_{k=1}^2\left(X^k\frac{\partial P'^{ij}}{\partial x^k}-P'^{kj}\frac{\partial X^i}{\partial x^k}-P'^{ik}\frac{\partial X^j }{\partial x^k} \right)=0\,,\qquad i,j=1,2,3,4
\]
Substituting into this equations  the following polynomial anzats
\begin{equation}\label{subs}
 P'^{ij}=\sum_{k,m=1}^2 a^{ij}_{km}(q)p_kp_m+\sum_{k=1}^2 b^{ij}_{k}(q)p_k+c^{ij}(q)\,.
\end{equation}
 for the entries of the bivector 
 \[
 P'=\sum P'^{ij} \frac{\partial}{\partial x_i}\wedge\frac{\partial}{\partial x_j}\,,
 \]
 we obtain  $60$ partial differential equations on 36 functions $a^{ij}_{km}(q),b^{ij}_{k}(q)$ and $c^{ij}(q)$ on two coordinates $q_1$ and $q_2$. 
  
These equations have been solved using various modern computer algebra systems for later comparison of the results obtained.
The final results obtained in this way were verified analytically.

\section{Non-integrable H\'{e}non-Heiles system}
The St\"{o}rmer-Verlet method applied to the Hamiltonian vector field (\ref{x-vec}) preserves a discrete version of the canonical symplectic form
\[
\omega=dp_1\wedge dq_1+dp_2\wedge dq_2\,,
\]
and, as sequence, the corresponding volume form $\Omega=\omega^2$ and the Poisson bivector $P=\omega^{-1}$. Combined with results from perturbation theory, this explains the excellent long-term behaviour of the method: long-term energy conservation, linear error growth and preservation of invariant tori in near-integrable systems, a discrete virial theorem and preservation of adiabatic invariants, see \cite{sv2003} and references therein. 

Below we introduce the second invariant symplectic form for the H\'{e}non-Heiles system, which could be useful to improve the properties of this and other numerical methods in application to this system. For generic potential $V(q_1,q_2)$ solution of the equations (\ref{g-inv})  depends on two parameters $a_{1,2}\in \mathbb R$
\[
P'=(a_1H+a_2)\,P\,.
\]
In some special cases, the solution of the equation  (\ref{g-inv})  depends on three or more parameters. As an example, for the   H\'{e}non-Heiles potential \cite{hh64}
\[
V(q_1,q_2)=q_1(a q_2^2 + b q_1^2)\,,
\]
with arbitrary $a$ and $b$, solution depends on three parameters 
\begin{equation}\label{p-hh}
P'=(a_1H+a_2)\,P+a_3\tilde{P}\,,
\end{equation}
where fourth rank skew symmetric tensor field $\tilde{P}$ of type (2,0) has the following form
\bq\label{pt-hh}
\tilde{P}=\left(
\begin{array}{cccc}
0&\frac{4 ( q_2  p_1- q_1  p_2)}{3}&  p_1^2- p_2^{2}-\frac{2}{3}  {a}  q_1   q_2^{2}+2  {b}   q_1^{3}
&2  p_1  p_2+\frac{8 a q_1^2  q_2 }{3}\\ 
\\
*&0&2  p_1  p_2+\frac{4}{3}  {a}   q_2^{3}+4  {b}   q_1^2  q_2&
 p_2^{2}- p_1^{2}+\frac{2}{3}  {a}  {q_1}   q_2^{2}-2  {b}   q_1^{3}\\ 
 \\
 *&*&0&
4  {a}  q_1  q_2p_1-\left(2  {a}  q_2^{2}+6  {b}  q_1^{2}\right)  p_2\\ \\
*&*&*&0
\end{array}\right)\,.
\eq
Multiplying invariant bivector $\tilde{P}$ (\ref{pt-hh}) on the scalar invariant $H^{2/3}$ we obtain Poisson bivector $\hat{P}=H^{2/3}\tilde{P}$ which satisfies to the Jacobi identity
\bq\label{jac-c}
[\![\hat{P},\hat{P}]\!]=0\,.
\eq
Here  $[\![.,.]\!]$ is a Schouten-Nijenhuis bracket  between  alternating multivector fields. The Poisson bivector $\hat{P}$ is also compatible with the bivector $H^{10/3}P$ obtained from the canonical Poisson bivector (\ref{p-can})
\[
[\![\hat{P},H^{10/3}P]\!]=0\,.
\]
The Poisson bivector $\hat{P}$ defines invariant symplectic form $\hat{\omega}=\hat{P}^{-1}=H^{-8/3}\tilde{w}$, where
\begin{align*}
\tilde{\omega}=&\left(\frac{aq_2^2 + 3bq_1^2}{2}p_2 - aq_1q_2p_1\right)dq_1\wedge dq_2+
\left(\frac{aq_1q_2^2}{6} - \frac{bq_1^3}{2} - \frac{p_1^2-p_2^2}{4}\right)dq_1\wedge dp_1\\
-&\left(\frac{aq_2^3}{3} + bq_1^2q_2 + \frac{p_1p_2}{2}\right)dq_1\wedge dp_2-\left(
\frac{2aq_1^2q_2}{3}+\frac{p_1p_2}{2}\right)dq_2\wedge dp_1\\
+&\left(-\frac{aq_1q_2^2}{6} + \frac{bq_1^3}{2} + \frac{p_1^2 - p_2^2}{4}\right)dq_2\wedge dp_2
+\frac{q_1p_2-q_2p_1}{3}dp_1\wedge dp_2\,.
\end{align*}
As a result we have two invariant symplectic forms $\omega$ and $\hat{\omega}$
\[
\mathcal L_X{\omega}=0 \qquad\mbox{and}\qquad
\mathcal L_X\hat{\omega}=0\,,
\]
so that standard volume form is equal to
\[
\Omega=\omega^2=4H^{10/3}\hat{\omega}^2\,.
\]
We can try to construct multisymplectic integrator for the H\'{e}non-Heiles system preserving both symplectic forms simultaneously. 

Bivector $\tilde{P}$ (\ref{pt-hh}) is the solution of equation $\mathcal L_X\tilde{P}=0$. Other properties of this bivector  $\tilde{P}$ are not pertinent to the construction of structure-preserving integrators. This bivector may be a Poisson bivector, which is incompatible with the canonical bivector $P$, or it may be a non-Poisson bivector et al. So, we consider only bivector $\tilde{P}$ and pre-symplectic form $\tilde{\omega}$ which are a well-defined polynomials on the whole phase space.
\begin{prop}
Invariant bivector $\tilde{P}$ (\ref{pt-hh}) satisfies equation
\bq\label{x2-eq}
\tilde{P}dH=2HX\,,
\eq
It means that vector field $Y=2HX$ has a bi-hamiltonian form
\[
Y=\tilde{P}dH=Pd\tilde{H}\,,
\]
where $P$ is canonical invariant bivector (\ref{p-can}) and $\tilde{H}=H^2$.
\end{prop}
The proof consists of a simple calculation.

Transformation of vector field $X$ to vector field $Y=2HX$ is associated to the change of time   
\[
dt\to 2Hd\tau
\]
for dynamical system  (\ref{m-eq}), and, therefore, we can say that tensor invariants $P$ and $\tilde{P}$ are "dependent" in a broad sense
\bq\label{pp-h}
(\tilde{P}-2HP)dH=0\,.
\eq
Similar change of time can be applied to other Hamiltonian systems. Nevertheless, it is easy to prove that
invariance equation (\ref{g-inv}) for the Hamiltonian vector field $X$  (\ref{x-vec}) with the following non weight-homogeneous polynomial potentials  
\[
V(q_1,q_2)=aq_1^3+bq_1q_2\,,
\]
has  solution depending only on two free parameters 
\[
P'=(a_1H+a_2)P\,,
\]
in contrast with the H\'{e}non-Heiles system.

That vector fields $Y=F(H)X$, where $F$ is a function of the Hamiltonian $H$, can be multi-Hamiltonian is quite natural fact. The fact that such vector fields are multi-Hamiltonian only for weight-homogeneous potentials is a more interesting. We can assume that this is caused by restricting the search space of solutions of the invariant equation (\ref{subs}) and restriction $F(H)=H$.

\section{First family  Hamiltonian flows preserving two symplectic forms}
If we put  $a=b=0$ in (\ref{pt-hh}) we obtain invariant bivector
\bq\label{pt-f}
P'_{h}=\left(\begin{array}{cccc}
0& \alpha (p_1q_2 - p_2q_1)& p_1^2 - p_2^2& 2p_1p_2\\
-\alpha(p_1q_2-p_2q_2)& 0& 2p_1p_2& p_2^2-p_1^2\\
p_2^2 - p_2^2& -2p_1p_2& 0& 0\\
-2p_1p_2& p_1^2 - p_2^2&0& 0
\end{array}\right)\,,\qquad \alpha\in\mathbb R\,,
\eq
for the free motion on the plane. The corresponding invariant 2-form looks like
\[
\omega'_h=(p_1^2-p_2^2)\left(dq_1\wedge dp_1-dq_2\wedge dp_2\right)+2p_1p_2\left(dq_1\wedge dp_2+dq_2\wedge dp_1\right)
-\alpha(p_1q_2-p_2q_1)dp_1\wedge dp_2\,.
\]
Additive deformations of this invariant bivector $P'_h$ associated with integrable Hamiltonian systems are discussed in \cite{ts11,ts11a,ts12}. Now we want to consider similar additive deformations for the non-integrable Hamiltonian systems.

\begin{prop}
If  potential $V(q_1,q_2)$ in the vector field $X$ (\ref{x-vec}) is equal to 
\[
V(q_1,q_2)=q_1^{4/\alpha}f\left(\frac{q_2}{q_1}\right)\,,\qquad \alpha\in\mathbb R\,,
\]
then the corresponding flow preserves canonical Poisson bivector $P$ (\ref{p-can}) 
and the following rank four bivector
\[
\tilde{P}=\left(
                  \begin{array}{cccc}
                    0 & \alpha(p_1q_2 - p_2q_1) & p_1^2-p_2^2-\alpha  q_2 \partial_2V+2 V &2p_1p_2+ \alpha q_1\partial_2 V \\ 
                    0 & 0 & 2p_1p_2+\alpha q_2\partial_1 V & p_2^2-p_1^2+\alpha  q_2 \partial_2V-2 V \\ 
                    * & * & 0 &2p_1\partial_2 V - 2p_2\partial_1V  \\ 
                   *  & * & * & 0 
                  \end{array}
                \right)
\]
which is an additive deformation of the Poisson bivector $P'_h$ (\ref{pt-f}). 
\end{prop}
In \cite{ts11,ts11a,ts12} the bivector $P_f$ was referred to as a geodesic bivector, whereas  additional part $\tilde{P}-P_f$ depending on $V$ was referred to as a potential bivector. 

For the proof we fix an entry of the bivector $\tilde{P}$.
\[\tilde{P}^{12}=\alpha(p_1q_2 - p_2q_1),\] 
while other entries of $\tilde{P}$ remain second order polynomials in momenta with coefficients depending on $q_{1,2}$ (\ref{subs}). At $\alpha\neq 0$ solutions of the invariance equation $\mathcal L_X\tilde{P}=0$ are equal to
\begin{align*}
\tilde{P}^{12}=& \alpha(p_1q_2 - p_2q_1)\,,\qquad  \tilde{P}^{24}=-\tilde{P}^{13}\,,\\ \\
\tilde{P}^{13}=&p_1^2-p_2^2-\alpha  q_2 \partial_2V(q_1,q_2)+2 V(q_1,q_2)+c_1 \,q_2^2+2 c_2 q_2+c_3\,,\\ \\
\tilde{P}^{14}=&2p_1p_2 + \alpha q_1\partial_2 V(q_1, q_2) -q_1 (c_1 q_2- c_2) + c_4q_2 + c_5\,,\\ \\
\tilde{P}^{23}=& 2p_1p_2 +\alpha q_2\partial_1 V(q_1, q_2) -q_1 (c_1q_2 - c_2) + c_4q_2 + c_5\,,\\ \\
\tilde{P}^{34}=&
\left(2\partial_2 V(q_1, q_2) + c_1q_2 + c_2\right)p_1 - \left(2\partial_1V(q_1, q_2) + c_1q_1 - c_4\right)p_2\,,
\end{align*}
where 
\begin{align*}
\alpha=2\,,\qquad
&V(q_1,q_2)=q_1^2f\left(\frac{q_2}{q_1}\right)
+c_1\left(\frac{(q_1^2+q_2^2)}{2}\ln q_1-\frac{q_1^2}{4}\right)- (c_2q_2-c_4q_1)\,,
\\
\alpha=4\,,\qquad
&V(q_1,q_2)=q_1f\left(\frac{q_2}{q_1}\right)+\frac{c_1(q_1^2+q_2^2)}{4}+\frac{c_2q_2-c_4q_1}{2}\,\ln q_1+\frac{c_2}{2}q_1
\\
\alpha\neq2\,,\alpha\neq 4\,,\qquad
&V(q_1,q_2)=q_1^{4/\alpha}f\left(\frac{q_2}{q_1}\right) + \frac{c_1(q_1^2 + q_2^2)}{2(\alpha - 2)}
 + \frac{2(c_2q_2 - c_4q_1)}{\alpha - 4}\,.
\end{align*}
Here constants of integration $c_1,\ldots,c_5$ and function $f(q_2/q_1)$ satisfy a total of one algebraic
\[c_1(c_4w + c_2)=0\]
and three differential equations, which  are omitted for the sake of brevity. 
In the event that all the constants of integration are set to zero, i.e. $c_1 = 0, \ldots, c_5 = 0$, then $f(q_2/q_1)$ is an arbitrary function. This provides a proof of the Proposition 2.

For the weight-homogeneous potential 
\[
V(q_1,q_2)=q_1^{4/\alpha}f\left(\frac{q_2}{q_1}\right)\,,\qquad \alpha\in\mathbb R\,,
\]
the invariance equation $\mathcal L_XP'=0$ (\ref{g-inv}) has the following solution depending on three parameters  $a_1,a_2$ and $a_3$
\[
P'=(a_1H+a_2)\,P+a_3\tilde{P}\,,\qquad a_k\in\mathbb R\,,
\]
Solving the corresponding Jacobi equation  $[\![P',P']\!]=0$ with respect to parameters
$a_k$ we obtain two nontrivial invariant Poisson bivectors $P'_{1,2}$ at
\[a_1 = 2a_3\,,\quad a_2=0\quad\mbox{and}\quad a_1 = \alpha a_3, \quad a_2=0.\]
In both cases we have the counterparts of equation (\ref{pp-h}) 
\bq\label{pp-hhg}
(P'-4a_3HP)dH=0 \qquad\mbox{and}\qquad
\Bigl(P' - a_3(\alpha + 2)HP\Bigr)dH=0\,,
\eq
respectively. In both cases, we have a pair incompatible Poisson bivectors $P$ and $P'$, i.e. $[\![P',P]\!]\neq0$. 
It is important to note that operator $N=P'P^{-1}$ is an invariant of the flow associated with $X$ together with its spectral data. However, it is crucial to note that its Nijenhuis torsion does not equal zero.

\section{Superintegrable Kepler problem} 
For degenerate in the Kolmogorov sense systems \cite{kolm27} there are global first integrals  which are functions on the both action and  angle variables \cite{ts08}. The Kepler problem is one of the most fundamental problems in physics having such  integrals of motion and, therefore, we take it as an example for our mathematical experiment. 

Following Euler \cite{eul} we immediately move to consider orbit plane dynamics with Cartesian coordinates $q_{1,2}$ so that the Hamiltonian $H$ and the corresponding vector field $X$ are
\begin{align}\label{x-kepl}
&H=\frac{p_1^2+p_2^2}{2}-\frac{\kappa}{\sqrt{q_1^2+q_2^2}}\,,\\ &X= p_1\dfrac{\partial}{\partial q_1}+p_2\dfrac{\partial}{\partial q_2}-
\frac{\kappa}{(q_1^2 + q_2^2)^{3/2}}\left(q_1\dfrac{\partial}{\partial p_1}+q_2\frac{\partial}{\partial p_2}\right)\,.\nonumber
\end{align}
This Hamiltonian commutes with the two components of the Laplace-Runge-Lenz vector
\[
K_1=p_1(p_1q_2 - p_2q_1)-\frac{\kappa q_2}{\sqrt{q_1^2 + q_2^2}}\,,\qquad
K_2= p_2(p_1q_2 - p_2q_1)+\frac{\kappa q_1}{\sqrt{q_1^2 + q_2^2}} \]
and the component of the orbital angular momentum
\[
K_3=q_1p_2-q_2p_1\,.\]
According to Euler \cite{eul} the pair of first integrals $H$ and $K_{1}$ (or $K_2$) has a St\"{a}ckel form  in elliptic coordinates on the orbit plane and the existence of an additional independent first integral $K_3$ is a consequence of the Euler additional law on elliptic curve  \cite{ts20}.

According to Jacobi \cite{jac-book} the pair of first integrals $H$ and $K_3$ has a St\"{a}ckel form in polar coordinates and
components of the Laplace-Runge-Lenz vector are derived using the Euler-Jacobi method of the last multiplier. 

Action-angle variables can be computed using both elliptic and polar coordinates. In this paper, we will limit ourselves to the consideration of the action-angle variables obtained using polar coordinates.

\subsection{Action-angle variables and known invariant bivectors}
Let us pass to the polar coordinates 
\[
q_1=r\cos\varphi\,,\quad q_2=r\sin\varphi\,\]
and the corresponding momenta
\[
p_1 = p_r \cos\varphi-\frac{p_\varphi\sin\varphi}{r}\,,\quad
p_2 = p_r\sin\varphi +\frac{p_\varphi\cos\varphi}{r}
\]
in which first integrals $H$ and $K_3^2$ have the St\"{a}ckel form \cite{st}
\begin{align*}H&=S^{-1}_{11}\left(p_r^2+V_1(r)\right)+S^{-1}_{21}\left(p_\varphi^2+V_2(\varphi)\right)=\frac{1}{2}\left(p_r^2 + \frac{p_\varphi^2}{r^2}\right)
 - \frac{\kappa}{r}\,,\\
 \\
K_3^2&=S^{-1}_{12}\left(p_r^2+V_1(r)\right)+S^{-1}_{22}\left(p_\varphi^2+V_2(\varphi)\right)=p_\varphi^2\,,
\end{align*}
where
\[S=\left(
 \begin{array}{cc}
 2 & 0 \\
 -r^{-2} & 1 \\
 \end{array}
 \right)\,,\qquad V_1(r)=-{2\kappa}{r}\,,\qquad V_2(\varphi)=0\,.
\]
According to \cite{min06} there is a discretization scheme which preserves these St\"{a}ckel integrals.

Next, for $H = h< 0$ we introduce action-angle variables 
\begin{align*}
I_\varphi&=p_\varphi\,,\qquad I_r = \frac{\kappa}{\sqrt{-2H}} - p_\varphi\,,\qquad H =-\frac{\kappa^2}{2(I_r+I_\varphi)^2}\\
\\
\theta_r& = \arctan\left(\frac{rp_r\sqrt{ 2\kappa r-p_r^2r^2 - p_\varphi^2}rp_r}{p_r^2 r^2 - \kappa r + p_\varphi^2}\right) - 
\frac{p_r\sqrt{2\kappa r-p_r^2r^2 - p_\varphi^2)}}{\kappa}\,,\\
\\
\theta_\varphi&=\theta_r + \varphi - \arcsin\left(\frac{\kappa r - p_\varphi^2}{\sqrt{p_\varphi^4 + r(p_r^2 r - 2\kappa)p_\varphi^2 + \kappa^2 r^2}}\right)
\end{align*}
and so-called Delauney elements \cite{mor02}
\[
I_1=I_\varphi,\quad I_2= I_r + I_\varphi,\quad \theta_1= \theta_\varphi - \theta_r,\quad \theta_2= \theta_r\,.
\]
Substituting the Delauney elements into the Bogoyavlenskij construction \cite{bog96,bog97} of the invariant 2-forms
\begin{equation}\label{inv-w}
 \omega'=\sum _{j=1}^2 d\left(\frac{\partial B(J)}{\partial J_j}\right)\wedge d\theta_j-df_j(I)\wedge d I_j\,,
 \end{equation}
 where $B(J_1,J_n)$ and $f(I_1,I_n)$ are arbitrary smooth functions  and
\[J_i=\frac{\partial H}{\partial I_i}\,,\qquad i=1,2\,.\]
we obtain a continuum of the local invariant Poisson bivectors $P'=\omega'^{-1}$ which are compatible  or non-compatible with canonical Poisson  bivector $P$. The main problem is to find global counterparts of these tensor invariants.

As an example, we present  one   invariant  Poisson bivector
\[
P_I'=I_1\frac{\partial}{\partial \theta_1}\wedge\frac{\partial}{\partial I_1}+I_2^{-2}\frac{\partial}{\partial \theta_2}\wedge\frac{\partial}{\partial I_2}
\]
which is a single valued tensor field in  polar variables and  the corresponding momenta
\begin{align}\label{kepl-p1}
P_I'=&\frac{2H}{\kappa^2}\frac{\partial}{\partial r}\wedge\frac{\partial}{\partial p_r}-p_\varphi\frac{\partial}{\partial \varphi}\wedge\frac{\partial}{\partial p_\varphi} \\
-&\frac{p_\varphi(\kappa^2 p_\varphi + 2H)}{\kappa^2 r(\kappa^2+2H p_\varphi^2)}\left(
(p_\varphi^2-\kappa r)\frac{\partial}{\partial r}\wedge\frac{\partial}{\partial \varphi}
+\frac{p_\varphi^2 p_r}{r}\frac{\partial}{\partial \varphi}\wedge\frac{\partial}{\partial p_r}
\right)\,.\nonumber
\end{align}
Singularity in the Poisson bivector $P_I'$ has no a physical meaning in contrast with systems considered in \cite{rat24}. Note that this invariant bivector $P_I'$ satisfies equation
\[
(P'_I+FP)dH\quad\mbox{with}\quad F=2\kappa^{-2}H\,,
\]
which is similar to equations (\ref{pp-h}) and (\ref{pp-hhg}).

Using other action-angle variables 
\[
J_{1,2}=\frac{1}{2}(H\pm I_1)\,,\qquad \chi_{1,2}=\frac{\theta_2I_2^3}{\kappa^2}\pm \theta_1
\]
we can rewrite Hamiltonian $H$ (\ref{x-kepl}) in the Fernandes form \cite{fer94}
\[H=J_1+J_2\]
and define so-called bi-Hamiltonian structure for the Kepler problem  associated with the Poisson bivector
\[P'=\left(
     \begin{array}{cccc}
       0 & 0 & J_1 & 0 \\
       0 & 0 & 0 & J_2 \\
       -J_1 & 0 & 0 & 0 \\
       0 &-J_2 & 0 & 0 \\
     \end{array}
   \right)\,,\qquad H=\frac{1}{2}\mbox{tr} P'P^{-1}\,.
  \] 
Recall, a bi-Hamiltonian system is prescribed by specifying two Hamiltonian functions $H_1$ and $H_2$
satisfying
\[X = PdH_1 = P'dH_2,\]
with diagonalizable recursion operator $N=P'P^{-1}$, having functionally independent real eigenvalues \cite{fer94}.
  
In polar variables and momenta this invariant Poisson bivector $P'$ are multivalued function on the phase space. So, we can say that bi-Hamiltonian structure  for the Kepler problem  in the Fernandes sense exists only locally  \cite{gts15}, i.e. in the neighbourhood of the open toroidal domains defined by the Arnold-Liouville theorem \cite{bog98}.  

\subsection{Mathematical experiment}

The substitution of polynomial anzats for entries of $P'$ and the Kepler vector field $X$ into the invariance equation $L_XP'=0$ (\ref{g-inv}) results in a system of 60 partial differential equations on 36 functions of $q_1$ and $q_2$.

\begin{prop} For  the Kepler system (\ref{x-kepl}) an invariance equation (\ref{g-inv}) has the following generic solution 
\begin{equation}\label{k-sol}
P'=(a_1 X_1+a_2X_2)\wedge X_3 + (a_3H+a_4K_3^2 + a_5K_1 +a_6K_2 + a_7K_3 + a_8)P+a_9\tilde{P}\,,
\end{equation}
depending on nine free parameters $a_i\in \mathbb R$  in the space of type (\ref{subs}) bivectors.  

Here $K_{1,2}$ are components of the Laplace-Runge-Lenz vector, $K_3$ is a component of the angular momentum vector, $X_k$ are the corresponding invariant vector fields
\[
X_1=PdK_1\,,\qquad X_2=PdK_2\,,\qquad X_3=PdK_3\,,
\] 
and entries of the  supplemental invariant bivector $\tilde{P}$ are equal to 
\begin{align*}
 \tilde{P}^{12}& =q_1p_2-p_1q_2\,,\quad \tilde{P}^{13}=-\frac{p_2^2}{2}+\frac{\kappa q_2^2}{(q_1^2 + q_2^2)^{3/2}}\,,\quad 
\tilde{P}^{14}=\frac{p_1p_2}{2}-\frac{\kappa q_1q_2}{(q_1^2 + q_2^2)^{3/2}}\,,\\
\\
 \tilde{P}^{23}& =\frac{p_1p_2}{2}-\frac{\kappa q_1q_2}{(q_1^2 + q_2^2)^{3/2}}\,,\quad
\tilde{P}^{24}=\frac{p_1^2}{2}+\frac{\kappa q_1^2}{(q_1^2 + q_2^2)^{3/2}}\,,\quad
\tilde{P}^{34}=\frac{\kappa(p_1q_2 - p_2q_1)}{2(q_1^2 + q_2^2)^{3/2}}\,.
\end{align*}
\end{prop}
Proof consists of a straightforward calculation.

Similar to the Hamilton-Jacobi equation  we can say that equation 
\[
\mathcal L_XP'=0
\]
has  a "complete integral" $P'$ (\ref{k-sol}) depending on a sufficient number of arbitrary constants,
 that allows to get all the integrals of motion from invariant (1,1) tensor field $N=P'P^{-1}$  
\[
\mbox{tr}\,N=2(a_9-2a_6)H +2 (a_2 - 2a_3)K_1 -2 (a_1+a_4)K_2 - 4a_5K_3^2 - 4a_7K_3 - 4a_8\,.
\]
In contrast with the standard recursion operator with vanishing Nijenhuis torsion  we obtain both commuting and non-commuting first integrals. It will be interesting to find similar generic invariant structure for a discrete scheme from \cite{kepl18}.

Substituting $P'$ (\ref{k-sol}) into the Jacobi identity 
\[[\![P',P']\!]=0\]
and solving the resulting equation for parameters $a_k$ we obtain four invariant  Poisson bivectors $P'_1,\ldots,P'_4,$ 
  \begin{align*} 
  P'_1&=a(X_1+iX_2)\wedge X_3+bP\,,\quad &\mbox{rank}P'_1=4\,,\\
  P'_2&=a(HP-2\tilde{P})\,,\quad &\mbox{rank}P'_2=4\,,\\
  P'_3&=a(X_1\wedge X_3-K_2P)+b(HP+\tilde{P})\,,\quad &\mbox{rank}P'_3=2\,,  \\
  P'_4&=a(X_2\wedge X_3+K_1P)+b(HP+\tilde{P})\,,\quad &\mbox{rank}P'_4=2\,,
  \end{align*}
so that
\begin{align*}
P'_1dH - \left(b - a(K_1 - iK_2)\right) PdK_3=&0\,,\qquad
2P'_3dH-bK_3PdH=0\,,\\
P'_2dH - K_3PdH - 3HPdK_3=&0\,,
\qquad 2P'_4dH-bK_3PdH=0\,.
\end{align*} 
Only one of these Poisson bivectors  is compatible with the canonical one
\[  [\![P,P_1']\!]=0\,,\qquad [\![P,P_k']\!]\neq 0\,,\quad k=2,3,4.\]
In polar coordinates this Poisson bivector are equal to
\begin{equation}\label{kepl-p2}
P'_1={\scriptstyle e^{i\varphi}}\left(
                            \begin{matrix}
                              0 & -p_\varphi & 0 & 0 \\
                              p_\varphi & 0 & \frac{i p_\varphi^2}{r^2} & -ip_r p_\varphi - \kappa +\frac{ p_\varphi^2}{r} \\
                              0 & -\frac{i p_\varphi^2}{r^2} & 0 & 0 \\
                              0 & ip_r p_\varphi + \kappa -\frac{ p_\varphi^2}{r} & 0 & 0 \\
                            \end{matrix}
                          \right)+bP
\end{equation}
In term of the Delauney elements it looks like
\[
P'_1=\phi_1 d\theta_1\wedge dI_1+ \phi_2 d\theta_1\wedge d\theta_2+bP
\]
where $\phi_{1,2}$ are functions on the Delauney elements which we  ommit for brevity. The  term involving angle variables
$d\theta_1\wedge d\theta_2$ is missing both in the Bogoyavlenski construction (\ref{w-inv}) in  \cite{bog96,bog97,bog98}, and in the bi-Hamiltonian geometry \cite{mv92}.

So, for the Kepler problem we have three invariant symplectic forms
\[
\omega=P^{-1}\,,\qquad \omega'_1={P'_1}^{-1}\,,\qquad \omega'_2={P'_2}^{-1}
\]
and two invariant rank-two Poisson bivectors $P_3$ and $P_4$.

\section{Second family of Hamiltonian flows preserving two symplectic forms}
Following to \cite{ts11,ts11a,ts12} we can substitute polynomial bivector $P'$ (\ref{subs}) into the invariance equation $\mathcal L_X P'=0$ and  solve a subsystem of the partial differential equations that do not include the potential $V(q_1,q_2)$ and its derivatives. As a result we obtain that following entry 
\[
P'^{12}=\sum_{k,m=1}^2 a^{12}_{km}(q_1,q_2)p_kp_m+\sum_{k=1}^2 b^{12}_{k}(q_1,q_2)p_k+c^{12}(q_1,q_2)
\]  
have to be 
\begin{align}
P'^{12}=&b_1 K^2 + (b_2 p_1 + b_3 p_2 + b_4 q_1 + b_5 q_2) K + b_6 p_1^2 +b_7 p_1 p_2 + b_8 p_2^2 + b_9 q_1^2 + b_{10} q_1 q_2 + b_{11} q_2^2\nonumber \\ \label{p12-gen}\\
+& (b_{12} q_1 + b_{13} q_2 + b_{14}) p_1 + (b_{15} q_1 + b_{16} q_2 + b_{17}) p_2 + q_1 b_{18} + q_2 b_{19} + b_{20}\,,
\nonumber
\end{align}
where $b_k\in\mathbb R$ and $K=(q_1p_2-q_2p_1)$. 

In the previous sections we study potentials $V(q_1,q_2)$ associated with the partial solution (\ref{pt-f}) of invariance equation \[P_h'^{12}=\alpha(q_1p_2-q_2p_1).\]
Now we wanted to prove that there are non-trivial potentials $V (q_1, q_2)$ associated with other simple partial solution of the invariance equation
\bq\label{pt-t}
P_t'^{12}=\alpha p_1+\beta p_1\,.
\eq
Following to \cite{ts11,ts11a,ts12}  we are looking for an invariant bivector  $\tilde{P}$ which is a sum of a geodesic bivector $P_t'$ defined by entry (\ref{pt-t}) and some additional bivector depending on potential $V$. 
\begin{prop}
If the Hamiltonian vector field $X$ (\ref{x-vec}) is defined by potential 
\bq\label{x-toda}
V(q_1,q_2)=\exp\left(\frac{-4q_1}{\beta}\right)f\left(\frac{\alpha q_1 + \beta q_2}{\beta}\right)\quad{or}\quad
V(q_1,q_2)=c_1\exp\left(\frac{-4q_1}{\beta}\right)+c_2\exp\left(\frac{4q_2}{\alpha}\right), 
\eq
then the corresponding flow preserves the canonical Poisson bivector $P$ (\ref{p-can}) 
and  the following rank four bivector 
\bq\label{pt-toda}
\tilde{P}=\left(
                  \begin{array}{cccc}
                    0 & \alpha p_1 + \beta p_2 & p_1^2 - p_2^2 - \alpha\partial_2 V + 2V & 2p_1p_2 - \beta \partial_2 V \\ 
                    0 & 0 & 2p_1p_2 + \alpha \partial_1 V & -p_1^2 + p_2^2 + \alpha\partial_2 V) - 2V) \\ 
                    * & * & 0 &2p_1\partial_2 V p_1 - 2p_2\partial_1V   \\ 
                   *  & * & * & 0 
                  \end{array}
                \right)\,.
\eq
Here $\alpha,\beta\in\mathbb R$ and $c_1,c_2\in\mathbb R$  are arbitrary parameters and $f\left((\alpha q_1 + \beta q_2)/\beta\right)$ is an arbitrary function.
\end{prop}
The proof is provided by the direct solution of the invariance equation for bivector with fixed entry (\ref{pt-t})
whereas other entries of $\tilde{P}$ remain the polynomials of second order in momenta with coefficients depending on $q_{1,2}$ (\ref{subs}). The integration constants arising in the solution process were assumed to be equal to zero. 

Bivector $\tilde{P}$ (\ref{pt-toda}) satisfies to the equation
\[
(\tilde{P}- 2HP) dH=0
\]
similar to equation (\ref{x2-eq}) for the H\'{e}non-Heiles system. For the potential 
\[
V(q_1,q_2)=\exp\left(\frac{-4q_1}{\beta}\right)f\left(\frac{\alpha q_1 + \beta q_2}{\beta}\right)\qquad \alpha,\beta\in\mathbb R\,,
\]
the invariance equation $\mathcal L_XP'=0$ (\ref{g-inv}) has the following solution depending on three parameters  $a_1,a_2$ and $a_3$
\[
P'=(a_1H+a_2)\,P+a_3\tilde{P}\,,\qquad a_k\in\mathbb R\,,
\]
Solving the corresponding Jacobi equation  $[\![P',P']\!]=0$ with respect to parameters
$a_k$ we obtain  nontrivial invariant Poisson bivector 
\[P' = 2HP+\tilde{P}\]
so that 
\bq\label{pp-ttg}
(P'-4HP)dH=0 
\eq
where $P$ and $P'$ are two incompatible Poisson bivectors $P'$ which are solution of the invariance equation. The corresponding vector field
\[
Y=P'dH=PdH'\,,\qquad H'=2H^2
\]
has a  bi-hamiltonian form. As sequence, operator $N=P'P^{-1}$ is an invariant of the flow,  but its Nijenhuis torsion does not equal to zero. 

In the generic case, the Hamiltonian vector field $X$ (\ref{x2-eq}) with potential (\ref{x-toda}) is not integrable. However, there are partial functions $f$ which yield integrable vector fields, in a manner analogous to the H\'{e}non-Heiles system. 

\subsection{Open and periodic Toda lattices}
Let us consider  Hamiltonian describing open Toda lattice ($\alpha=0$) and periodic Toda lattice ($\alpha=1$) associated with the root system $G_2$
\bq\label{g-toda}
H=\frac{p_1^2 + p_2^2}{2} + \exp(q_1/\sqrt{3}) + \exp(-\sqrt{3}/2q_1 + q_2/2) + \alpha\exp(-q_2)\,.
\eq
A detailed description of these systems can be found in \cite{g2}. We only recall that the second integral of motion for both these integrable systems are polynomials of sixth order in momenta.

\begin{prop}
For the open Toda lattice with $\alpha=0$ in (\ref{g-toda}) the invariance equation $\mathcal L_XP'=0$ has the following solution
\begin{equation}\label{sol-toda}
P'=(a_1H+a_2)P+a_3\tilde{P}\,,\qquad a_k\in\mathbb R\,,
\end{equation} 
depending on three parameters.The entries of the  invariant bivector $\tilde{P}$ are equal to
\begin{align*}
\tilde{P}^{12}&=\sqrt{3}p_2 - 5p_1\,,\qquad\qquad \tilde{P}^{13}=\frac{p_2^2}{2} + \frac{5\exp(-\sqrt{3}/2q_1 + q_2/2)}{2}\,,
\\
\tilde{P}^{14}&=-\frac{p_1p_2}{2} - \frac{\sqrt{3}\exp(-\sqrt{3}/2q_1 + q_2/2)}{2}\,,
\\
\tilde{P}^{23}&=-\frac{p_1p_2}{2} -\frac{5\sqrt{3}\exp(q_1/\sqrt{3})}{3} +\frac{5\sqrt{3}\exp(-\sqrt{3}/2q_1 + q_2/2)}{2}\,,
\\
\tilde{P}^{24}&=\frac{p_1^2}{2} -\frac{ 3\exp(-\sqrt{3}/2q_1 + q_2/2)}{2} +\exp(q_1/\sqrt{3})\,,
\\
\tilde{P}^{34}&=
-\frac{\exp(-\sqrt{3}/2q_1 + q_2/2)}{4}p_1 \\
&+\frac{\sqrt{3}\left(2\exp(q_1/\sqrt{3}) - 3\exp(-\sqrt{3}/2q_1 + q_2/2)\right)}{12}p_2\,.
\end{align*}
For the periodic Toda lattice (\ref{g-toda}) the invariance equation $\mathcal L_XP'=0$ has solution
\[ P'=(a_1H+a_2)P \]
depending only on two parameters in the  space of bivectors with entries (\ref{subs}).
\end{prop}
Proof consists of a straightforward calculation.

The most interesting result of this experiment is the fact that  non-trivial solution  of the invariance equations
(\ref{g-inv}) exists for an open Toda lattice ($\alpha=0$) and not for a periodic lattice ($\alpha=1$).

\section{Conclusion}
The subject of geometric numerical integration is concerned with numerical integrators that preserve geometric properties of the flow of a differential equation \cite{mar2001}. For instance, each finite-dimensional Hamiltonian flow preserves the Hamiltonian vector field $X$ defined by
\[
\iota_X\omega=dH
\] 
together with the symplectic form $\omega$ and the Hamiltonian $H$, i.e.
\[\mathcal L_XX=0\,,\qquad \mathcal L_X\omega=0\,,\qquad \mathcal L_XH=0\,.\]
 The variational  energy-preserving  integrators ensure that the numerical solution is confined to a codimension one submanifold  of the configuration manifold defined by equation $H=E$. In contrast, variational symplectic integrators maintain the canonical symplectic form $\omega$ at each time step.
A comparison of the numerical performance of both methods for benchmark problems from various types of Hamiltonian systems can be found in  \cite{book0,book1,book2}. 

Multi-symplectic integrators are usually regarded as a discretisation that conserves a discrete version of the conservation of symplecticity for Hamiltonian partial differential equations (PDEs). This is because, for generic finite-dimensional Hamiltonian systems (ODEs), the construction of the multi-symplectic structure is based on the same basic invariants $\omega$ and $H$ \cite{leon99,bar18}. We are aware that there are at least two invariant symplectic forms only for integrable bi-Hamiltonian systems \cite{bog96,bog97,bog98}. 
 
In this note we present  second invariant symplectic form $\omega'$ for some well-studied Hamiltonian systems with two degrees of freedom: nonintegrable H\'{e}non-Heiles system, Kepler problem, integrable and non-integrable Toda type lattice
\[\mathcal L_X\omega'=0\,.\]
It allows us to construct two symplectic integrators preserving invariant symplectic forms $\omega$ and $\omega'$, respectively, and to understand which integrator performs better for the given Hamiltonian systems. Furthermore, we can attempt to construct a multi-symplectic integrator that simultaneously preserves both symplectic forms $\omega$ and $\omega'$.

In accordance with the Darboux theorem \cite{mbook}, there exist coordinates $q',p'$ in which the second invariant symplectic form $\omega'$ assumes its canonical form
\[
\omega=dp_1'\wedge dq_1'+dp'_2\wedge dq'_2\,.
\] 
It is of interest to study property of this transformation $(q,p)$ to $(q',p')$ for nonintegrable Hamiltonian systems.

For Hamiltonian systems on the Poisson manifolds,  solutions of the invariance  equation  (\ref{g-inv}) 
are more diverse. In particular, these solutions generate  polynomial brackets on the Lie algebras of small dimension and we can try to classify them similar to the Cartan classification of the linear Lie brackets \cite{ts25}.

The article was prepared within the framework of the project “International academic cooperation” HSE University.

\end{document}